\documentclass[12pt]{article}
\usepackage{amsmath, amssymb, amsfonts, amsthm, amscd, exscale, enumerate}
\usepackage[latin1]{inputenc}

\newcommand{\RR}{\mathbb{R}}

\newtheorem{theorem}{Theorem}
\newtheorem{condition}[theorem]{Condition}
\newtheorem{corollary}[theorem]{Corollary}
\newtheorem{definition}[theorem]{Definition}

\newtheorem{lemma}[theorem]{Lemma}

\newtheorem{proposition}[theorem]{Proposition}
\newtheorem{remark}[theorem]{Remark}

\begin{document}

\title{A Hardy inequality on Riemannian manifolds and a classification of discrete Dirichlet spectra}
\bigskip
\author{Nils Rautenberg\\Ruhr-Universit\"at Bochum, Fakult\"at f\"ur
Mathematik,\\ Universit\"atsstr. 150, D-44780 Bochum, Germany.\\
Email : Nils.Rautenberg@ruhr-uni-bochum.de} \maketitle

\begin{abstract} 
\noindent We prove a Hardy inequality for uniformly elliptic operators subject to Dirichlet or mixed boundary conditions on domains $\Omega$ with piecewiese smooth boundary in arbitrary Riemannian Manifolds $(M,g)$. Employing an approach of E.B. Davies for the euclidean case, we show that it implies a sufficient geometric criterion under which the Laplace-Beltrami operator with Dirichlet boundary conditions $\Delta^D$ has purely discrete spectrum on $\Omega$. We proceed to classify all non-compact $\Omega$ with discrete spectrum up to a boundary regularity condition and show that these include for example polygons in manifolds of negative curvature with ideal vertices. This a new result for non-constant curvature.\end{abstract}
$\\$

\newpage
\section{Introduction and statement of the results}
Let $(M^n,g)$ be a complete smooth Riemannian manifold of dimension $n$. Let $\Omega \subset M$ be a domain, that is, an open connected subset with piecewise smooth boundary. More precisely, the boundary $\partial \Omega$ is the disjoint union of a (possibly infinite) number of $(n-1)$-dimensional smooth open pieces and its complement, the set of critical points $C$, which is assumed to have measure zero in $\partial \Omega$.

 Let $H$ denote a uniformly elliptic second order operator acting on $L^2(\Omega,d\mu)$ with potential bounded from below subject to Dirichlet boundary conditions. So $H=A+V$ with 
\[ \int\limits_\Omega fA(f) d\mu\geq \alpha \int\limits_\Omega f \Delta (f) d\mu. \] 
For some $\alpha>0$. Here $d\mu$ denotes the usual Riemannian measure induced by $g$. Let $Q$ denote the quadratic form of $H$, so in particular for all $f \in \mathrm{Dom}(H)$ we have:

\[Q(f)=\int\limits_\Omega f H(f) d\mu . \]
Let $\delta$ be a differentiable function on $\Omega$ satisfying $||\mathrm{grad}(\delta)||_g\leq 1$. $H$ satisifes a (weak) Hardy inequality with respect to $\delta$ if there exist $a\geq 0, c>0$ such that
\[ c\left( Q(f) + a\int\limits_\Omega |f|^2 d\mu\right)  \geq  \int\limits_\Omega \frac{|f|^2}{\delta^2}d\mu\]
for all $f \in H^1_0:=\overline{C_0^\infty(\Omega)}$, where the closure is taken with respect to the Sobolev norm $\|\cdot\|_1$ given by $\|f\|^2_1=\int\limits_\Omega g(\mathrm{grad}(f),\mathrm{grad}(f)) d\mu + \int\limits_\Omega |f|^2d\mu$.\\

If $a=0$, then such an inequality is called a strong Hardy inequality and $c$ is then called a strong Hardy constant. Many versions of such inequalities have been established in different contexts. A review can be found for example in \cite{Dav2}. \\
One inequality of this kind (see for example \cite{Dav}) is the following: Let $\Omega \subset \RR^n$ be a domain and let $H=\Delta:=-\mathrm{div}~\mathrm{grad}$ be the Laplacian. Then for all $f \in H^1_0(\Omega)$ the following strong Hardy inequality holds:

\begin{equation} Q(f)=\int\limits_\Omega||\mathrm{grad}(f)||^2  d\mu \geq \frac{n}{4} \int\limits_\Omega \frac{|f|^2}{m^2}d\mu,
\end{equation}
where the weight $m(x)$ is defined by
\[\frac{1}{m(x)^2}=\int\limits_{S^{n-1}}\frac{1}{r_x(v)^2}dv , \quad \quad \mathrm{with} \quad  r_x(v):=\inf \{ |t|~:~ x+tv \in \partial \Omega\}\]
and $dv$ denotes the normalized euclidean surface measure of the sphere. That is, $m(x)$ represents a mean distance to the boundary evaluated along all lines through $x$. By definition, $m(x)\geq d(x):=d(x,\partial \Omega)$ holds.\\

The main result of this paper is that a generalized version of inequality (1) holds true on any Riemannian manifold by replacing lines with geodesics and by replacing the average over the sphere by an average over the unit tangent space $S_p\Omega$ at each point $p \in \Omega$. 

\begin{theorem}
Let $\Omega \subset (M^n,g)$ be a domain with piecewise smooth boundary in a complete Riemannian manifold of dimension $n$. Then, for all $f \in H_0^1(\Omega)$:
\[\int\limits_\Omega||\mathrm{grad}(f)||_g^2  d\mu \geq \frac{n}{4} \int\limits_\Omega \frac{|f|^2}{m^2}d\mu\]
and
\[\frac{1}{m(p)^2}=\int\limits_{S_p\Omega}\frac{1}{r_p(v)^2}dv, \quad \quad \mathrm{with} \quad r_p(v):=\inf \{ |t|~:~ c_v(t) \in \partial \Omega\},\]
where $c_v$ is the geodesic with initial conditions $c_v(0)=p$ and $\dot{c}_v(0)=v$ and $dv$ denotes the normalized euclidean surface measure on $S_p\Omega\cong S^{n-1}$.
\end{theorem}
Note that the constant does not depend on the sectional curvature $K$ or other geometric data and is explicit as a function of the dimension only. However, the function $m$ can usually not be computed directly. Assume that $H$ is uniformly elliptic with potential bounded from below, then a direct consequence of Theorem 1 is the following:

\begin{theorem}
 If $H$ is uniformly elliptic with potential bounded from below, then for all $f \in H^1_0(\Omega)$ the (weak) Hardy inequality:
\[Q(f)+a\int\limits_\Omega|f|^2d\mu\geq \frac{\alpha n}{4} \int\limits_\Omega \frac{|f|^2}{m^2}d\mu\]
holds with $a=\max\{0, -\inf(V(x)) \}$ and $\alpha$ the constant of uniform ellipticity of $H$.
\end{theorem}

As in the euclidean case, a large class of domains permits an estimate of the form $d(p)\leq m(p) \leq c d(p)$ with $c>1$ for all points $p\in \Omega$, where $d(p):=d_g(p,\partial \Omega)$ is now the Riemannian distance to the boundary. We call such domains \emph{boundary distance regular}.\\

Let $\sigma(H)$ denote the spectrum of the self-adjoint operator $H$. We have the standard decomposition $\sigma(H)=\sigma_{disc}(H)\cup~\sigma_{ess}(H)$. The discrete spectrum $\sigma_{disc}(H) $ consists of all isolated eigenvalues of finite multiplicity, and the essential spectrum is the closure of the complement of $\sigma_{disc}(H)$ in $\sigma(H)$.\\

If one considers the Laplace-Beltrami operator $\Delta=- \mathrm{div_g}~\mathrm{grad_g}$ with Dirichlet conditions, denoted by $\Delta^D$, then Theorem 1 implies a classification of all boundary distance regular domains such that on them the spectrum of $\Delta^D$ is purely discrete, that is, $\sigma_{ess}(\Delta^D)=\emptyset$. A domain $\Omega \subset (M^n,g)$ is said to be $\emph{quasi-bounded}$ iff for every $\epsilon>0$ the set $\Omega_\epsilon:=\{p \in \Omega~|~ d(p)\geq\epsilon \}$ is compact. We prove the following theorem for $\Delta^D$:
\begin{theorem}
Let $\Omega \subset (M^n,g)$ be a boundary distance regular domain. Then the quasi-boundedness of $\Omega$ is a sufficient condition for purely discrete spcetrum of $\Delta^D$, in other words $\sigma_{ess}(\Delta^D) =\emptyset$.
\end{theorem}
\begin{remark}
 If $K>-a$ for some $a>0$, then quasi-boundedness of $\Omega$ is necessary for purely discrete spectrum of $\Delta^D$.
\end{remark}
This is well-known and follows for example from min-max arguments and Cheng's theorem \cite{Che} that allows one to compare eigenvalues on geodesic balls to eigenvalues of geodesic balls on spaces of constant curvature. \\

Boundary distance regularity has been investigated thoroughly in the past, see for example \cite{Dav}, p. 27. We state two conditions that imply it, though this is not an exhaustive list:

\begin{condition} (Uniform Interior Cone (UIC) Condition)\\
Assume that for $\Omega \subset (M,g)$ there exists an angle $\alpha>0$ and a constant $c_0>1$ such that for each $p \in \Omega$ there exists an $\alpha$-angled cone $C_{\alpha}(p)\subset T_p\Omega$ of directions with the property that $r_p(w)< c_0 d(p)$ for each $w \in C_{\alpha}$. Then $\Omega$ is boundary distance regular.
\end{condition}
To see this, one estimates the integral in the definition of $m$ from below by an integral of $c_0d(p)$ over the cone.

\begin{condition} (Uniform Exterior Ball Condition (UEB) Condition)\\
Suppose that there exists a constant $k>0$ such that for each $p \in \partial \Omega$, there exists a ball $B(q,kd_g(p,q))$ around a point q and of radius $kd_g(p,q)$ disjoint from $\Omega$, then $\Omega$ is boundary-distance regular.
\end{condition}
See for example Theorem. 1.5.4 in \cite{Dav}. \\

In the penultimate section, we will proceed to apply this classification theorem to an example in the form of polygons $P$ in negatively curved simply connected manifolds, which will have at least one vertex on the ideal boundary. We will show that they are both boundary distance regular and quasi-bounded and thus have purely discrete spectrum of $\Delta^D$. This was not known outside constant negative curvature.\\

In the final section of this paper we will establish a version of theorem 1 that applies to operators with mixed boundary conditions and gives a similar inequality for their quadratic forms.\\

We conclude the introduction by remarking that our method of proof is inspired by the work of Croke and Derdzinski \cite{CD}. Croke proved that for every complete Riemannian manifold with boundary $(M^n,g)$ the following bound holds for the first eigenvalue $\lambda_1(M,g)$ of $\Delta^D$:
\begin{equation} \lambda_1(M,g) \geq \frac{n \pi}{\mathrm{vol}(S^{n-1})} \inf\limits_{p \in M} \int\limits_{SpM}\frac{1}{l(v)^2}dv
\end{equation}
where $l(v)$ is the length of the maximal geodesic $c_v$. 
Furthermore, it should be noted that the inequality of Croke is already sharp. Specifically, he and Derdzinski proved in the same paper that equality is equivalent to $M$ being a Riemannian hemisphere bundle.\\

Though unlinke inequality (2), the Hardy inequality of Theorem 1 allows no direct evaluation of the integral even in special cases, it does allow for better qualitative statements, since for $v \in S\Omega$, it may be that $l(v)$ is very large, even infinite, while $r_p(v)$ is small.

\section{Preliminaries}

Let $(M^n,g)$ denote a complete smooth Riemannian manifold of dimension $n$. We define $T_pM$ to be the tangent space at a point $p \in M$, and $S_pM$ to be the space of tangent vectors at $p$ of length one.\\

Furthermore, with $I$ being an intervall containing $0$, define by $c_v: I \to M$ the geodesic with intital conditions $c_v(0)=p$, $\dot{c}_v(0)=v$, with $v \in S_pM$. We also consider $\pi:TM \to M$ resp. $\pi:SM\to M$, the tangent bundle and unit tangent bundle of $M$. We will write an element $(p,v) \in SM$ simply as $v$, making the convention that for $v \in SM$, $p:=\pi(v)$. On $M$, the geodesic flow is defined as the smooth one-parameter group of diffeomorphisms given by:

\begin{eqnarray*}
\phi^t: SM &\to& SM \\
            v &\mapsto& \phi^t(v):=\dot{c}_v(t)
\end{eqnarray*}

Consider a domain $\Omega \subset M$ and its closure as a metric space $\bar{\Omega}$. If the boundary $\partial \Omega$ is not empty, a geodesic in $\Omega$ might hit this boundary in finite time and in such a way that it does not exist beyond this point as a geodesic in $\Omega$. But we can always think of it as a geodesic segment in $M$, possessing a natural and unique extension to a geodesic in $M$ that exists for all times. We will not differ between the two in our notation, as the maximal geodesic in $\Omega$ is merely a restriction of the maximal geodesic in $M$.\\

Assume that the boundary is piecewise smooth. Denote the complement of the smooth pieces by $C$. Throughout this paper, we will assume that $C$ is a set of measure zero within $\partial \Omega$. Let $N$ denote the inward pointing unit normal vector field on the smooth boundary components. Define:

\[S^+\partial \Omega:=\{u \in SM ~|~ \pi(u)=p \in \partial \Omega \setminus C ~\mathrm{and}~ g_p(N,u)>0\}\] and
\[S^-\partial \Omega:=\{u \in SM ~|~ \pi(u)=p \in \partial \Omega \setminus C ~\mathrm{and}~ g_p(N,u)<0\}\] 
Both sets are $(2n-2)$-dimensional submanifolds of $SM$. Consider $u \in S^+\partial \Omega$. Then the maximal geodesic $c_u$  in $\Omega$ is defined on an interval $I_u$ which is either of the form $(0, \infty)$ or $(0,l(u))$, with $0<l(u)<\infty$. Note that in the former case, we must have $c_u(l(u)) \in \partial \Omega$ and $g_{c_u(l(u))}(N, \dot{c}_u(l(u)))\leq 0$ whenever $c_u(l(u)) \notin C$. Note also that if $l(u)=\infty$, then $c_{-u}$ will be defined on $(-\infty, 0)$ as a geodesic in $\Omega$. The same holds, up to different signs, for $u \in S^-\partial \Omega$. 
Let $A:=\bigcup\limits_{u \in S^+\partial \Omega}\{u\}\times I_u$, then the geodesic flow gives a map:

\begin{eqnarray*}
F:~~ A &\to& S\Omega \\
            (u,t) &\mapsto& F(u,t):=\phi^t(u)=\dot{c}_u(t)
\end{eqnarray*}

If we consider the set $S^*\Omega:=\{ v \in S\Omega ~|~ \exists t_0> 0 ~\mathrm{s.t.}~ r_p(v) < \infty\}$, then the image $F(A)$ will be an open subset of positive measure in $S^*\Omega$. Furthermore, the map is injective and smooth on $A$. Let $d\mu_L$ denote the volume form associated with the canonical Liouville measure on $SM$ as well as its restriction onto $S^*\Omega$ and the image $F(A)$. We have the following formula for the pullback of the volume form $F^*d\mu_l$:

\begin{theorem} (Santalo's formula)\\
Let $du$ be the volume form on $S^+\partial \Omega$ induced by the volume form on $SM$. Then we have:
 \[ F^*d\mu_L=g(N,u)dt \wedge du, \]
where $dt$ is the canonical measure on $\mathbb{R}$.
\end{theorem}
\paragraph{Proof:}

See \cite{San}, p. 337 or \cite{Ber} p. 282-285. $\hfill \square$

Note that the formula also holds true if we define the map $F$ on the vectors of $S^-\partial \Omega$ and their intervals of existence. Next, let us define the set:
\[ S^{-}_\infty\partial \Omega:=\{ u \in SM ~|~ -u \in S^+\partial \Omega ~\mathrm{and}~ I_{-u}=(0, \infty) \}\]
and finally  consider the union 
\[S^*\partial \Omega:=S^+\partial \Omega \cup S^-_\infty \partial \Omega\]

The map $F$ extends to $S^*\partial \Omega$ and is injective. The image of $F$ on this larger set is dense in $S^*\Omega$ since we exclude only the critical points $C$ on the boundary and the directions tangential to the boundary. Since $S^-_\infty \partial \Omega$ is a measurable subset of $S^-\partial \Omega$, we can restrict the measure induced by the pull-back of the volume form $F^*d\mu_L$ to this subset and the Santalo formula still holds true as a formula for the induced measure on $S^*\Omega$.\\

In order to finish the preparation of the proofs, we state the original Hardy inequality in dimension 1:
\begin{lemma} (Hardy's inequality)\\
 If $f:[a,b] \to \mathbb{C}$ is continuously differentiable with $f(a)=0=f(b)$, then:
\[\int\limits_a^b\frac{|f(x)|^2}{4d(x)^2}dx \leq \int\limits_a^b |f^\prime(x)|^2dx\]
with $d(x)= \mathrm{min} \{|x-a|,|x-b|\}$
\end{lemma}
We give the proof, as stated in \cite{Dav}, p. 26:
\paragraph{Proof:}
It is sufficient to prove:
\[ \int\limits_a^c\frac{|f(x)|^2}{4(x-a)^2}dx \leq \int_a^c|f^\prime(x)|^2dx,\]
where $2c=a+b$, and a similar inequality for the other half-interval. It is sufficient to deal with the case where $a=0$, and where $f$ is real. Then:
\begin{eqnarray*}
\int\limits_0^c(f^\prime)^2dx&=&\int\limits_0^c\left(x^\frac{1}{2}(x^{-\frac{1}{2}}f)^\prime+\frac{f}{2x}\right)^2 dx\geq \int\limits_0^c\left( x^{-\frac{1}{2}}f(x^{-\frac{1}{2}}f)^\prime+\frac{f^2}{4x^2}\right) dx \\
&=& \left[\frac{1}{2}(x^{-\frac{1}{2}} f)^2\right]_0^c+ \int\limits_0^c\frac{f^2}{4x^2}dx=\frac{f(c)^2}{2c}+\int\limits_0^c\frac{f^2}{4x^2}dx\geq \int\limits_0^c\frac{f^2}{4x^2}dx
\end{eqnarray*}
$\hfill \square$\\
Note that the proof immediatly yields the following:
\begin{corollary} (Hardy's inequality for half-axes)\\
 If $f:[a,\infty) \to \mathbb{C}$ is continuously differentiable and square integrable with $f(a)=0$, then:
\[\int\limits_a^\infty\frac{|f(x)|^2}{4|x-a|^2}dx \leq \int\limits_a^\infty |f^\prime(x)|^2dx.\]
\end{corollary}

\section{Proof of the main theorems:}
\paragraph{Proof of Theorem 1:} 
To begin with, note that for all $f \in C^\infty_0(\Omega)$ we have the equality:
\[ \|\mathrm{grad} f(p)\|_g^2=\frac{n}{\mathrm{vol}(S^{n-1})}\int\limits_{Sp\Omega} (df(p)(v))^2 dv \]
where $dv$ is the canonical spherical suface measure. Using this, we find:

\begin{eqnarray*}
\int\limits_\Omega \|\mathrm{grad} f(p)\|_g^2d\mu(p)&=&\frac{n}{\mathrm{vol}(S^{n-1})}\int\limits_\Omega\int\limits_{Sp\Omega}df(p)(v)^2dvd\mu(p)\\
&=&\frac{n}{\mathrm{vol}(S^{n-1})}\int\limits_{S\Omega}df(\pi(v))(v)^2d\mu_L(v)\\
&\geq&\frac{n}{\mathrm{vol}(S^{n-1})}\int\limits_{S\Omega^*}df(\pi(v))(v)^2d\mu_L(v),\\
\end{eqnarray*}

where we remind that $S^*\Omega$ denotes all $v \in S\Omega$ for which $r_p(v)<\infty$. Using Santalo's formula and Hardy's inequality, we derive for $f \in C^\infty_0(\Omega)$:

\begin{eqnarray*}
&&\frac{n}{\mathrm{vol}(S^{n-1})}\int\limits_{S^*\Omega}df(\pi(v))(v)^2d\mu_L(v)\\
&=&\frac{n}{\mathrm{vol}(S^{n-1})}\int\limits_{S^*\partial \Omega}\left(\int\limits_{0}^{l(u)}df(\pi(\phi^t(u)))(\phi^t(u))^2 dt \right) g_{\pi(u)}(N(\pi(u)),u) du\\
&=&\frac{n}{\mathrm{vol}(S^{n-1})}\int\limits_{S^*\partial \Omega}\left(\int\limits_{0}^{l(u)}d(f \circ \pi(\phi^t(u)))^2 dt \right) g_{\pi(u)}(N(\pi(u)),u) du\\
&\geq&\frac{n}{\mathrm{vol}(S^{n-1})}\int\limits_{S^*\partial \Omega}\left( \int\limits_{0}^{l(u)}\frac{f(\pi(\phi^t(u)))^2}{4r_{\pi(\phi^t(u))}(\phi^t(u))}dt\right) g_{\pi(u)}(N(\pi(u)),u) du\\
&=&\frac{n}{\mathrm{vol}(S^{n-1})}\int\limits_{S^*\Omega}\frac{f(\pi(v))^2}{4r_p(v)^2}d\mu_L(v)
\end{eqnarray*}
\begin{eqnarray*}
&=&\frac{n}{\mathrm{vol}(S^{n-1})}\int\limits_{S\Omega}\frac{f(\pi(v))^2}{4r_p(v)^2}d\mu_L(v)\\
&=&\frac{n}{4}\int\limits_\Omega\frac{f(p)^2}{m(p)^2} d\mu(p),
\end{eqnarray*}
where the second-last equality holds since $r_p(v)=\infty$ for all $v \in S\Omega\setminus S^*\Omega$. Since $C_0^\infty(\Omega)$ is dense in $H^1_0(\Omega)$ we are done by Fatou's lemma. $\hfill \square$

\paragraph{Proof of Theorem 2:}
As a direct consequence of the definition of $H=A+V$, we know

\[H\geq\alpha \Delta^D+V\]
in the sense of quadratic forms. If the potential is non-negative, one omits it and as a result derives a strong Hardy inequality, while in the case where $V$ is negative but bounded, one just adds $a=- \inf (V(p))$ to the quadratic form and finds a weak Hardy inequality. $\hfill \square$
\paragraph{Proof of Theorem 3:}
Assuming boundary distance regularity of $\Omega$,  the Hardy inequality implies that one has the following inequality in the sense of quadratic forms for any $b>0$:
\[ \Delta^D\geq \frac{1}{2} \left(\Delta^D+ \frac{nc}{4d^2}\right)=\frac{1}{2}\left(\Delta^D+V_b+W_b\right)\]
With 
\begin{equation*}
V_b(p)=\begin{cases} b, &\text{if}~ \frac{nc}{4d(p) ^2}\leq b \\ \frac{nc}{4d(p)^2}, & \text{else} \end{cases} 
\end{equation*}
 and 
\begin{equation*}
W_b(p)=\begin{cases}  \frac{nc}{4d(p)^2}-b, &\text{if}~ \frac{nc}{4d(p) ^2}\leq b \\ 0, &\text{else}.
\end{cases} 
\end{equation*}

This is an equality between quadratic forms, and, by for example min-max arguments, we can use to bound the ssential spectrum of $\Delta^D$ from below by the right hand side. Since $\Omega$ is quasi-bounded, $W_b$ is continous with compact support. Then, by standard arguments the difference of the resolvents of $\Delta^D+V_b+W_b$ and $\Delta^D+V_b$ is a compact operator. A theorem of Weyl (cf. \cite{RS}, Theorem XIII.14) then yields:

\[ \sigma_{ess}(\Delta^D+V_b+W_b)=\sigma_{ess}(\Delta^D+V_b) \]

However, the spectrum of the latter operator is contained in $[b,\infty)$, since the Laplacian is positve and the potential is bounded from below by $b$. Because $b$ was arbitrary, the essential spectrum of $\Delta^D$ is empty as claimed. $\hfill \square$

\section{Application: $\sigma(\Delta^D)$ on hyperbolic polygons with ideal vertices}

Let $<\cdot,\cdot>$ denote the euclidean metric, $||\cdot||$ its induced norm and $D^{n-1}$ the open euclidean unit ball. Then $(M,g)=(D^{n-1}, \frac{4}{((1-||p||^2)^2}<\cdot,\cdot>)=:\mathbb{H}^n$ is the Poincare ball model of hyperbolic $n$-space.\\
We consider a k-sided geodesic polygon $P\subset \mathbb{H}^2$, for which at least one of the interior angles is zero, or equivalently, that has at least one vertex on the ideal boundary $\partial \mathbb{H}^2\cong S^1$. Though it is not a bounded set in the hyperbolic metric, $P$ is of finite volume. We show the following well-known statement with the help of Theorem 3:

\begin{proposition}
The spectrum of $\Delta^D$ on any poylgon $P\subset \mathbb{H}^2$ with ideal vertices is discrete.
\end{proposition}

\paragraph{Proof:}
We have to verify that $P$ is (i) quasi-bounded and (ii) boundary distance regular. 

\begin{itemize}
\item[i)] For quasi-boundedness we observe it suffices to show that for any $\epsilon>0$ no ideal vertex is in the closure of $P_\epsilon$. On each ideal vertex, the 2 bounding geodesics forming the vertex are asymptotic to each other in the hyperbolic metric. Therefore for each $\epsilon>0$ the ideal vertex cannot be a limit point of the set of points with a distance to both geodesics of at least $\epsilon$, of which $P_\epsilon$ is a subset. Thus, all $P_\epsilon$ are compact as claimed.

\item[ii)]
To show that $P$ is boundary distance regular, we verify that they fullfill the uniform interior cone condition (UIC). Note that the distance of any point in the interior of the polygon to the boundary is uniformly bounded and that the maximum $d_0$ is achieved at at least one point $p_0$. For any $p$ in the interior of the polygon the boundary distance $d(p)$ is realised by a geodesic arc $c_v$ that hits the boundary at a right angle. 
We will prove the existance of a cone $C_\alpha(p)$ of directions around each such $v$, such that for each $w \in C_\alpha(p)$ the geodesic $c_w$ will hit the boundary after a time of at most $2d(p)$ and then bound $\alpha$ away from $0$ uniformly for all $p \in P$.\\

Since the boundary $\partial \Omega$ consists of geodesic arcs in $\mathbb{H}^2$, the existence of $C_\alpha(w) $ follows from the existence of certain right angled geodesic triangles. One side is given by a portion of one of the boundary geodesics, another side is given by the distance realizing geodesic arc $c_v$ and the hypothenuse is given by a geodesic arc such that the ratio of the length of the hypothenuse and length of $c_v$ is a number $1<r\leq 2$. Of course, the choice of a factor of $2$ was arbitrary and the existance of such triangles for all $r$ is a well-known geometric fact. The inital conditions for the hypothenuses then form $C_\alpha(p)$.\\

It remains to argue why the angle of these cones can be bounded uniformly from below. This, however follows from the geometric fact that for a fixed ratio of  lengths $r$, the angle between the hypothenuse and the arc $c_v$ increases as the length decreases. Since the distance to the boundary is bounded from above by $d_0$, we set the uniform cone angle to be the angle of one of the cones around a distance realising arc $c_{v_0}$ at $p_0$. Note that at $p_0$, there might be more than one direction belonging to such a distance realising arc. However, for each of these the cone will be of equal size by the homogenity of $\mathbb{H}^2$. Thus, (UIC) holds. 
\end{itemize}
 $\hfill \square$

In the same way, one can argue in higher dimensions, where a $k$-sided polygon $P \subset \mathbb{H}^n$ is now given as a set that is bounded by $k$ totally geodesic subspaces $\mathbb{H}^{n-1}$, of which at least two are asymptotic to each other at the ideal boundary. The argument remains unchanged, and we note:\\

\begin{corollary}
The spectrum of $\Delta^D$ on a polygon $P \in \mathbb{H}^n$ with ideal vertices is discrete.
\end{corollary}

Furthermore, we can generalize the argument to variable negative curvature. That is, consider a simply connected manifold $M$ of dimension $n$ that satisfies the sectional curvature bounds $-b\leq K \leq -a$ with $b>a>0$ and define an $k$-sided polygon $P$ with ideal vertex to be a subset of $M$ bounded by $k$ totally geodesic submanifolds $N_i$, $1\leq i \leq n$, of dimension $n-1$, of which at least two are asymptotic at the ideal boundary.\\

\begin{proposition}
On any simply connected Riemanannian manifold $(M,g)$ with curvature $K$ satisfying $-b\leq K \leq -a$, $b>a>0$, the spectrum of $\Delta^D$ on a Polygon $P$ is discrete even if it has ideal vertices.
\end{proposition}

\paragraph{Proof:}
Using the triangle comparison theorem of Topogonov \cite{Pet}, compare right angled geodesic triangles in $P$ with those of the constant curvature space form $M_{-b}$. Both quasi-boundedness and boundary distance regularity follow directly from the properties of the polygons in $M_{-b}$.
$\hfill \square$

\section{Extension to operators with mixed boundary conditions}
In this section we show how to extend Theorem 1 and 2 to operators with mixed boundary conditions. The weight $m$ is then replaced by a restriction to an average over those directions that 'see' the Dirichlet component of the boundary. More precisely:

\begin{definition}
Let $\Omega \subset (M,g)$ be a domain in a Riemannian manifold. Let $\Gamma \subset \partial \Omega$ be a closed subset of the boundary. Define $m_\Gamma : \Omega \to (0,\infty) \cup \{\infty\}$ by:
\[\frac{1}{m_\Gamma^2(p)}:=\int\limits_{S_p\Omega}\frac{1}{r^\Gamma_p(v)^2}dv\]
with $r^\Gamma_p(v):= \inf \{ |t|~:~c_v(t) \in \Gamma \}.$
\end{definition}
If we furthermore define $H^1_{0,\Gamma}(\Omega):=\{ f\in H^1(\Omega) ~:~ f_{|\Gamma} \equiv 0\}$, we have the following generalization of Theorem 1:

\begin{theorem}Let $\Omega \subset (M^n,g)$ be a domain in a Riemannian manifold. Let $\Gamma \subset \partial \Omega$ be a piecewise smooth closed subset of the boundary. For all $f \in H^1_{0,\Gamma}(\Omega)$ the following inequality holds:
\[\int\limits_\Omega||\mathrm{grad}(f)||_g^2  d\mu \geq \frac{n}{4} \int\limits_\Omega \frac{|f|^2}{m_\Gamma^2}d\mu\]
\end{theorem}
\paragraph{Proof:} The proof of Theorem 1 holds line by line using $m_\Gamma$ instead of $m$, and working in $H^1_{0,\Gamma}(\Omega)$ instead of $H^1_0(\Omega)$. $\hfill \square$\\

Since $H^1_{0,\Gamma}$ is the form domain of $\Delta^{DN}$, the Laplacian with Dirichlet boundary conditions on $\Gamma$ and Neumann conditions everywhere else,  we have completed our task. It is clear from the definition of $H$ and the inequality for the mixed Laplacian that Theorem 2 generalizes in just the same way to include uniformly elliptic operators with mixed boundary conditions.

\end{document}